\documentclass[12pt]{article}

\usepackage{graphicx,amsmath,amssymb,subfig}
\usepackage[affil-it]{authblk}

\begin{document}

\title{Clustering-based collocation for uncertainty propagation with multivariate dependent inputs}
\author[1]{A.W. Eggels \thanks{a.w.eggels@cwi.nl}}
\author[1,2]{D.T. Crommelin}
\author[1]{J.A.S. Witteveen}
\affil[1]{Centrum Wiskunde \& Informatica, Amsterdam, the Netherlands}
\affil[2]{Korteweg - de Vries Institute for Mathematics, University of Amsterdam, the Netherlands}


	\maketitle
	\noindent Original submission: 29 March 2017 \\ Revised on: 9 October 2017 \\ Author's original manuscript of \\ DOI:10.1615/Int.J.UncertaintyQuantification.2018020215

\abstract{In this article, we propose the use of partitioning and clustering methods as an alternative to Gaussian quadrature for stochastic collocation. The key idea is to use cluster centers as the nodes for collocation. In this way, we can extend the use of collocation methods to uncertainty propagation with multivariate, dependent input, in which the output approximation is piecewise constant on the clusters. The approach is particularly useful in situations where the probability distribution of the input is unknown, and only a sample from the input distribution is available.
	We examine several clustering methods and assess the convergence of collocation based on these methods both theoretically and numerically. We demonstrate good performance of the proposed methods, most notably for the challenging case of nonlinearly dependent inputs in higher dimensions. Numerical tests with input dimension up to $16$ are included, using as benchmarks the Genz test functions and a test case from computational fluid dynamics (lid-driven cavity flow). 
}
\section{Introduction}
A core topic in the field of uncertainty quantification (UQ) is the question how to characterize the distribution of model outputs, given the distribution of the model inputs (or a sample thereof). Questions such as these are encountered  in many fields of science and engineering  \cite{Bij13,Wal02,Wit08,Wit07,Yil15}, and have given rise to modern UQ methods including stochastic collocation, polynomial chaos expansion and stochastic Galerkin methods \cite{Xiu02,Xiu05,Gha03,Mai10,Eld09}.

A still outstanding challenge is how to characterize model output distributions efficiently in case of multivariate, dependent input distributions. In the previously mentioned methods independence between the inputs is assumed, e.g., for the construction of the Lagrange polynomials in stochastic collocation, or for the construction of the orthogonal polynomials in generalized polynomial chaos. When independence between the input components holds, the multivariate problem can easily be factored into multiple $1$-dimensional problems, whose solutions can be combined by tensor products to a solution for the multidimensional problem. When the inputs are dependent, such factorization can become extremely complicated if the inputs have non-Gaussian distributions, making it unfeasible in practice for many cases. It generally involves nontrivial transformations that require detailed knowledge of the joint distribution (e.g. Rosenblatt transformation \cite{2Ros52}), however such information is often not available. In \cite{Nav15}, factorization is circumvented and instead the problem is tackled by using the Gram-Schmidt (GS) orthogonalization procedure to get an orthogonal basis of polynomials, in which the orthogonality is with respect to the distribution of the inputs. However, this procedure gives non-unique results that depend on the implementation.

In this paper we propose a novel approach for efficient UQ with multivariate, dependent inputs. This approach is related to stochastic collocation, however it employs collocation nodes that are obtained from data clustering rather than from constructing a standard (e.g. Gaussian) quadrature or cubature rule. By using techniques from data clustering, we can construct sets of nodes that give a good representation of the input data distribution, well capable of capturing correlations and nonlinear structures in the input distributions. It is straightforward to obtain weights associated with these nodes. All weights are guaranteed to be positive.

The approach we propose is non-intrusive and able to handle non-Gaussian dependent inputs. We demonstrate that it remains efficient for higher dimensions of the inputs, notably in case of strong dependencies. These dependencies are not limited to correlations (linear dependencies), but can also be nonlinear. Furthermore, the approach employs data clustering, starting from a sample dataset of inputs. The underlying input distribution can be unknown, and there is no fitting of the distribution involved. Thus, no fitting error is introduced. This makes the approach particularly suitable for situations where the exact input distribution is unknown and only a sample of it is available.

We emphasize that the method we propose in this paper does not employ orthogonal polynomials and their roots, nor does it require to specify an input distribution. This constitutes a main difference from stochastic collocation. Furthermore, we demonstrate that generating a random quadrature rule, by randomly selecting points from the sample of inputs and using these as cluster centers, gives unsatisfactory results. This is due to the fact that such a random selection is ill-suited to sample or represent the tails of the input distribution.

The outline of this paper is the following: in Section \ref{sec:scext}, we start by briefly summarizing stochastic collocation and multivariate inputs. We discuss the challenges of dealing with dependent inputs, and we introduce the concept of clustering-based collocation. In Section \ref{sec:methods}, we describe three different clustering techniques and give a convergence result for one dimension. In Section \ref{sec:res}, we present results of numerical experiments in which we test our clustering-based collocation method, using the clustering techniques described in Section \ref{sec:methods}. A test case from computational fluid dynamics (lid-driven cavity flow) is described in Section \ref{sec:ldcf}. The conclusion follows in Section \ref{sec:conclusion}.

\section{Stochastic collocation and its extension}\label{sec:scext}
Consider a function $u(\mathbf{x}): \Omega\mapsto\mathbb{R} \, , \,\,\,\, \Omega\subseteq\mathbb{R}^p$, that maps a vector of input variables to a scalar output. Let us assume $\mathbf{x}$ is a realization of a random variable $\chi$ with probability density function $f(\mathbf{x})$. We would like to characterize the probability distribution of $u(\mathbf{x})$, in particular we would like to compute moments of $u(\mathbf{x})$:
\begin{equation}
\mathbb{E}[u^q] = \int_\Omega (u(\mathbf{x}))^q f(\mathbf{x})d\mathbf{x} \, .
\end{equation}
In what follows, we focus on the first moment:
\begin{equation}
\label{eq:mean}
\mu := \mathbb{E} \, u = \int_\Omega u(\mathbf{x}) f(\mathbf{x})d\mathbf{x} \, .
\end{equation}
We note that higher moments can be treated in the same way, as these are effectively averages of different output functions, i.e. $\mathbb{E} [u^q] = \mathbb{E} v$ with $v(\mathbf{x}) := (u(\mathbf{x}))^q$. In both cases, the expectation is with respect to the distribution of $\chi$. In stochastic collocation, the integral in (\ref{eq:mean}) is approximated using a quadrature or cubature rule. As is well-known, a high degree of exactness of the integration can be achieved for polynomial integrands with Gaussian quadrature rules.

\subsection{Multivariate inputs}
For multivariate inputs ($p>1$), stochastic collocation based on Gaussian quadrature can be constructed using tensor products if the input variables are mutually independent. In this case, we can write $f(\mathbf{x})$ as a product of $1$-dimensional probability density functions. The degree of exactness of the corresponding cubature rule is $2k-1$ in each dimension if $k$ collocation nodes are used in each input dimension. This requires a number of nodes ($k^p$) that grows exponentially in $p$, so collocation with full tensor grids suffers from curse of dimension. To reduce the number of nodes, Smolyak sparse grids \cite{Smo63,2Ger98} can be used. The construction of Smolyak sparse grids will not be explained in detail here, but an important aspect is that the resulting set of nodes is a union of subsets of full tensor grids. When the grids are nested and the number of nodes in the $n$th level for one dimension, $k_n^1$, is $O(2^n)$, then the number of nodes in $p$ dimensions scales as $O(2^n n^{p-1})$ \cite{2Ger98}. This in contrast to $O(2^{np})$ for the corresponding tensor rule.

\subsection{Gaussian cubature with dependent inputs}
As already mentioned, tensor grids are useful for stochastic collocation in case of independent inputs. If the input variables are dependent, grids constructed as tensor products of $1$-dimensional Gaussian quadrature nodes no longer give rise to a Gaussian cubature rule. In \cite{Nav15}, generalization to dependent inputs is approached by constructing sets of polynomials that are orthogonal with respect to general multivariate input distributions, using Gram-Schmidt orthogonalization. The roots of such a set of polynomials can serve as nodes for a Gaussian cubature rule.

With the approach pursued in \cite{Nav15}, the advantages of Gaussian quadrature (in particular, its high degree of exactness) carry over to the multivariate, dependent case. However, one encounters several difficulties with this approach. First of all, for a given input distribution, the set of nodes that is obtained is not unique. Rather, the resulting set depends on the precise ordering of the monomials that enter the GS procedure. For example, with $2$-dimensional inputs and cubic monomials, $24$ different sets of nodes can be constructed, as demonstrated in \cite{Nav15}. It is not obvious a priori which of these sets is optimal.

A further challenge is the computation of the weights for the cubature rule. It is not straightforward how to construct multivariate Lagrange interpolating functions and evaluate their integrals. The alternative for computing the weights is to solve the moment equations. However, the resulting weights can be negative. Furthermore, one cannot choose the number of nodes freely: in general, with input dimension $p$ and polynomials of degree $m$, one obtains $n=m^p$ nodes. Thus, the number of nodes increases in large steps, for example with $p=8$ the number of nodes jumps from $1$ to $256$ to $6561$, respectively, if $m$ increases from $1$ to $2$ to $3$. It is unknown how to construct useful (sparse) subsets of nodes from these.

\subsection{Clustering-based collocation}
To circumvent the difficulties of Gaussian cubature in case of dependent inputs, as summarized in the previous section, we propose an alternative approach to choose collocation nodes. By no longer requiring the collocation nodes to be the nodes of an appropriate Gaussian cubature rule, we do not benefit anymore from the maximal degree of exact integration associated with Gaussian quadrature or cubature. However, we argue below that this benefit of Gaussian cubature offers only limited advantage in practice.

If one has a sample of the inputs available but the underlying input distribution is unknown, the Gaussian cubature rule will be affected by the sampling error (via the GS orthogonalization). Alternatively, if the input distribution is estimated from input sample data, the precision of the Gaussian cubature rule is also limited by the finite sample size.

Additionally, the degree of exactness is strongly limited by the number of nodes in higher dimensions. For example, suppose one can afford no more than $256$ evaluations of the output function $u(\mathbf{x})$ because of high computational cost, i.e. one can afford a Gaussian cubature rule with $256$ nodes. This gives very high degree of integration exactness (degree $511$) in one dimension ($p=1$), but the degree of exactness decreases to $31$, $7$ and $3$, respectively, as the input dimension $p$ increases to $2$, $4$ and $8$. For $p>8$, the degree of exactness is only $1$ in case of $256$ nodes, so only linear functions can be integrated exactly. The number of nodes for a full tensor grid in $p$ dimensions with $2^n$ nodes in level $n$ for one dimension is $2^{np}$ ($O(2^{np})$), while a corresponding Smolyak grid contains a number of points in the order $O(2^n n^{p-1})$. However, the approximation accuracy for the full grid is $O(2^{-nm})$ and $O(2^{-nm}n^{(p-1)(m+1)})$ for the sparse grid with $O(2^n)$ nodes in level $n$ for one dimension \cite{2Ger98}, where $m$ is the smoothness of the function. This is still limiting for a high number of dimensions.

Furthermore, the accuracy of the propagation method does not need to be higher than the accuracy of the input uncertainty. Since the input is given by samples, the high accuracy of spectral methods is not fully utilized.

Instead of constructing a Gaussian cubature rule, we aim to determine a set of nodes that are representative for the sample of input data or for the input distribution, with the locations of the nodes adjusting to the shape of the distribution. Clustering is a suitable method (or rather collection of methods) to achieve this objective. More specifically, clustering is the mathematical problem of grouping a set of objects (e.g., data points) in such a way that objects in one group (or cluster) are more similar to each other than to objects in other clusters \cite{Ste56,Jai99}. For each of the clusters, a center is defined to represent the cluster.

The basic idea, in the context of this study, is the following. Assume we have a dataset $\{ \mathbf{x}_1, ..., \mathbf{x}_N\}$ available, with $\mathbf{x}_i \in \mathbb{R}^p$. We define a partitioning of $\mathbb{R}^p$ existing of $K$ subsets, denoted $\Omega_k$ with $k=1, ..., K$. A cluster is a subset of the data falling into the same $\Omega_k$. A common way to define cluster centers $z_k$ is as the average of the data in each cluster, i.e.
\begin{equation}
\label{eq:def_zk}
\mathbf{z}_k := \frac{\sum_{i=1}^N \mathbf{x}_i \, \mathbf{1} (\mathbf{x}_i \in \Omega_k)}{\sum_{i=1}^N \mathbf{1} (\mathbf{x}_i \in \Omega_k)} \, ,
\end{equation}
with $\mathbf{1} (\cdot)$ the indicator function. If we define weights $w_k$ as the fraction of all the data falling in the $k$-th cluster, that is,
\begin{equation}
\label{eq:def_wk}
w_k : = \frac{\sum_{i=1}^N  \mathbf{1} (\mathbf{x}_i \in \Omega_k)}{N} \, ,
\end{equation}
the weighted average $\bar{\mathbf{z}} := \sum_k w_k \, \mathbf{z}_k$ equals the data average $\bar{\mathbf{x}} := N^{-1} \, \sum_i \mathbf{x}_i$. Thus, $\bar{\mathbf{z}}  = \bar{ \mathbf{x}}$  by construction.

The key idea of what we propose here is to carry out collocation based on clustering of the input data. More specifically, we propose to use the cluster centers $\mathbf{z}_k$ and weights $w_k$ as the nodes and weights of a quadrature rule. Thus, the (exact) first moment of the output function $u(\mathbf{x})$ over the input data is
\begin{equation}\label{eq:MC}
\mu = \frac1N \sum_{i=1}^N u(\mathbf{x}_i)
\end{equation}
and the approximation using clustering-based collocation is
\begin{equation}
\label{eq:quad_cbc}
\hat{\mu} := \sum_{k=1}^K w_k \, u(\mathbf{z}_k).
\end{equation}
We emphasize that the number of function evaluations in Equation (\ref{eq:MC}) and (\ref{eq:quad_cbc}) is different. When $K<<N$, large savings in computational time can be achieved due to the smaller amount of evaluations of $u(\mathbf{x})$.

The proposed approximation (\ref{eq:quad_cbc}) to estimate the first moment of $u(\mathbf{x})$ does not explicitly consider a function approximation of $u(\mathbf{x})$. However, (\ref{eq:quad_cbc}) can be seen as the Monte Carlo integral over a function approximation of $u(\mathbf{x})$ which is piecewise constant on the clusters.

It is easy to show that the approximation is exact ($\hat \mu = \mu$) for all linear input functions, due to the fact that $\bar{ \mathbf{z}}  = \bar{\mathbf{x}}$, as mentioned above. In other words, the degree of exactness is one: we can consider (\ref{eq:quad_cbc}) as a quadrature rule for the integral of $u(\mathbf{x})$ over the empirical measure induced by the dataset $\{\mathbf{x}_1, ..., \mathbf{x}_N\}$. This quadrature rule is exact if $u(\mathbf{x})$ is linear. This may seem limited in comparison to Gaussian quadrature, however as discussed earlier, the degree of exactness of Gaussian quadrature reduces rapidly if the input dimension $p$ grows and the number of nodes remains constant. For non-linear input functions, the approximation (\ref{eq:quad_cbc}) will in general not be exact. However, we will investigate its convergence in Sections \ref{ssec:convergencetheory} and \ref{sec:res}.

\section{Clustering methods}\label{sec:methods}
In this section, we describe three different methods to construct clusters, i.e. three methods to construct a suitable collection of subsets $\Omega_k$.
As already mentioned, the methods are based on input given as a dataset in $p$ dimensions with $N$ data points  $\{\mathbf{x}_1,\ldots,\mathbf{x}_N\}$ with $\mathbf{x}_i\in\mathbb{R}^{1\times p}$ also denoted by a matrix $X\in\mathbb{R}^{N\times p}$. If the input is given as a distribution, we can create a dataset by sampling from this distribution. Furthermore, we scale this dataset to $[0,1]^p$ by linear scaling with the range. This is done to comply with the domain of the test functions we will use further on.

We cluster the data points into $K$ clusters $\{C_1,\ldots,C_K\}$ with centers $\{\mathbf{z}_1, \ldots, \mathbf{z}_K\}$ and use these as nodes. The centers are computed as the mean of the data points in that cluster, see (\ref{eq:def_zk}). We investigate three different methods, namely $k$-means clustering, principal component analysis based clustering and a method with randomly selected data points as cluster centers. In the following, we will use the words clustering and partitioning interchangeably.

\subsection{$K$-means}\label{ssec:kme}
The $k$-means method is one of the oldest and most widely used methods for clustering \cite{Ste56,Mac67}. The idea behind it is to minimize the within-cluster-sum of squares (SOS):
\begin{equation}
\min_{\{\mathbf{z}_1,\ldots,\mathbf{z}_K\}} SOS(\mathbf{z}_1,\ldots,\mathbf{z}_K), \quad SOS(\mathbf{z}_1,\ldots,\mathbf{z}_K) = \sum_{i=1}^N ||\mathbf{x}_i-\mathbf{z}_{\text{argmin}_k ||\mathbf{x}_i-\mathbf{z}_k||_2^2}||_2^2.
\end{equation}
The minimization problem is solved with an iterative procedure, see e.g. \cite{Mac67} for details. There are many extensions and improvements of the (initialization of the) algorithm, such as using the triangle inequality to avoid unnecessary distance calculations \cite{Elk03}, the use of global methods \cite{Lik03,Bag08,Han05} and low-rank approximations \cite{Coh14}. We will use the $k$-means++ method in this subsection, which has a special initialization as described in \cite{Art07}.

Because the algorithm contains a random initialization and the objective function is non-convex, it can converge to a local minimum, rather than to the global optimum. Therefore, in our numerical tests in Section \ref{sec:res}, the algorithm is performed $r$ times $(r>1)$ with different initializations and the best solution (with minimal SOS) is chosen. We use a fixed number of iterations in the minimization. In some cases, the iterations have not fully converged yet.
This will be ignored because in practice, nearly all of the $r$ executions converge so that the chosen best solution is always a converged solution. Further onwards, we will refer to this method as KME. We choose $r=25$.

\subsection{PCA-based clustering}\label{ssec:pca}
With this method, based on \cite{Din04,Su04} and principal component analysis (PCA), one starts with a single large cluster containing all the data, and in each step, the cluster with the largest average radius is split in two. This is implemented by splitting the cluster whose data points have the largest average squared distance to the cluster center. We split such that the cutting plane goes through the old cluster center (center of mass) and is perpendicular to the largest principal component of the covariance matrix of the data in the cluster, as suggested by \cite{Su04}. This continues until the desired number of clusters $K$ is attained (we note that other stopping criteria can be used as well, however these are less useful for the purpose of this study).

This method (referred to as PCA later on) is deterministic, unlike the $k$-means method described in the previous section. Clustering by the diameter criterion is already performed in \cite{Gue91}, but there the cluster with the largest diameter is split. Division methods based on farthest centroids have been suggested by \cite{Fan08}. Other refinements of this method are also possible, e.g., the merging of clusters at some steps in the algorithm, but we will not explore these here. These can be investigated in future work.

\subsection{Random clustering}
For comparison purposes we include a third method, in which cluster centers are selected randomly. This method consists of randomly selecting data points from the data set, all with equal probability, and use these as cluster centers. The clusters are formed by assigning each data point to its nearest cluster center. This method will be referred to as MCC (Monte Carlo clustering).

\subsection{Calculation of weights}\label{ssec:weights}
As already mentioned, we use the cluster centers as nodes for collocation. To do so, each node must be assigned a weight.
In all three methods, the weight of each node is determined by the number of data points in the cluster associated with that node, divided by the total number of data points, see (\ref{eq:def_wk}). By construction, all weights are positive and their sum equals one.

\subsection{Convergence}\label{ssec:convergencetheory}
In the case of one dimension ($p=1$), it can be proven that the PCA-based clustering converges to the Monte Carlo integral for increasing values of $K$. The proof relies on the fact that in each step, the largest cluster radius either decreases or remains constant. If $p=1$, the largest cluster radius equals
\begin{equation}
\delta^*(K) = \max_{i\in\{1,\ldots,N\}}\min_{j\in\{1,\ldots,K\}}\{|x_i-z_j|\}.
\end{equation}
As can be seen, it depends on $K$. We give the proof for one dimension.

We can define a finite interval $D:=[x_{-},x_{+}]$ which contains all the data $\{ x_1, ..., x_N \}$. Furthermore, we assume that the output function (denoted $f(x)$ in this section) is Lipschitz continuous on this interval with Lipschitz constant $L$. Let $\nu$ be the empirical measure on $D$, i.\ e.\   $\nu(\Omega\subseteq D) = \frac{1}{N}\sum_{i=1}^N \mathbf{1}(x_i\in\Omega)$ for each subset $\Omega$ of $D$.
Suppose that we have for each $K\in\mathbb{N^+}$, $K\leq N$ that $x_{-}<z_1<z_2<\ldots<z_{K}<x_{+}$ are the ordered cluster centers. A set partitioning is defined by $\cup_{j=1}^K E_j$, where $E_j = \{x_i\in C_j\}$, for $j=1,\ldots,K$ with $C_j := \left[\frac{z_{j-1}+z_j}{2},\frac{z_{j}+z_{j+1}}{2}\right)$, for $j=2,\ldots,K-1$, $C_1=\left[x_{-},\frac{z_1+z_2}{2}\right)$ and $C_K=\left[\frac{z_{K-1}+z_K}{2},x_{+}\right]$.
Define $f_k(x):=f(z_k)$ $\forall x\in C_k$ and $0$ elsewhere for $k=1,\ldots,K$. Now, denote $\tilde{f}(x)=\sum_{k=1}^{K} f_k(x)$. Because of the Lipschitz continuity, we have that
\begin{equation}\label{eq:Lipschitz}
\forall x_i \in D \exists k =k(x_i)\in\{1,\ldots,K\}: |f(x_i)-\tilde{f}(x_i)| = |f(x_i)-f(z_{k})|<L\delta^*(K) 
\end{equation}
In the PCA-algorithm for $p=1$, $\delta^*$ is strictly non-increasing as $K$ grows. It reaches its lower bound $\delta^*(K)=0$ when $K=N$, because then each data point is its own cluster center.
We can now bound the difference between the PCA integral $I_{PCA}(K)=\sum_{k=1}^Kf(z_k)w_k$ and the Monte Carlo integral $I_{MC}(N)=\sum_{i=1}^N f(x_i)\frac{1}{N}$ as follows
\begin{align}
\left|I_{PCA}-I_{MC}\right| & = \left|\frac{1}{N}\sum_{i=1}^N f(x_i) - \sum_{k=1}^{K}f(z_k)w_k\right| \nonumber\\
& = \left|\frac{1}{N}\sum_{i=1}^N f(x_i) - \frac{1}{N}\sum_{i=1}^N\sum_{k=1}^{K}f(z_k)\mathbf{1}(x_i\in
\Omega_k)\right| \nonumber\\
& = \frac{1}{N}\left|\sum_{i=1}^N \left(f(x_i) - \sum_{k=1}^{K}f(z_k)\mathbf{1}(x_i\in
\Omega_k)\right)\right| \nonumber\\
& \leq \frac{1}{N}\sum_{i=1}^N \left|f(x_i) - \sum_{k=1}^{K}f(z_k)\mathbf{1}(x_i\in
\Omega_k)\right| \nonumber\\
& < \frac{1}{N}\sum_{i=1}^N L\delta^*(K) = L\delta^*(K).
\end{align}
Since $\delta^*(k+1)\leq \delta^*(k)$ for all $k\in\mathbb{N}$ when $p=1$, $\delta^*(k)\geq 0$ and $\delta^*(k=N)=0$, the bound becomes stricter for increasing $k$.

For higher dimensions, the derivation of the bounds is analogous, although $\delta^*$ will not be monotonically decreasing, but it will decrease in general. This is also the case for the other methods, even for $p=1$, where the nodes are not nested such as in the PCA-case, such that it is not guaranteed that $\delta^*$ decreases monotonically. For the higher-dimensional case, we want to refer to Section \ref{ssec:k} where we give a numerical result on convergence.

\newpage
\section{Results}\label{sec:res}
We test the quadrature based on the clustering methods described in Section \ref{sec:methods} by integrating the Genz test functions on the domain $[0,1]^p$ for three different data sets. In two of these data sets, the variables are mutually dependent. The relative error, as defined by the absolute difference between the integral calculated by the cluster points and weights and the Monte Carlo integral of the data, is used as the measure of accuracy.
We perform the MCC method $10$ times for each setting to investigate the effect of randomness. We show the mean, minimum and maximum error for MCC. For comparison, we have also added results from using Monte Carlo sampling (MCS), which was repeated $10$ times as well.

The first test is to assess how these methods perform under an increasing number of dimensions and what the effect of dependent variables is. Then, we compare the numerical convergence of the PCA method with the MCC and MCS method for an increasing number of clusters. Finally, we compare the computational cost of these methods.

\subsection{Genz test functions}
Genz \cite{Gen84} has developed several functions to test the accuracy of a cubature rule. The definitions, our choice of parameters and some illustrations are given in \ref{app:genz}. We test the methods by integrating the Genz test functions over the three different data sets consisting of $N=10^5$ samples. We compare the results from Monte Carlo integration and with the results obtained by clustering-based quadrature. The difference between the two integrals is a measure for the (in)accuracy of the methods.

\subsection{Data sets}
We use three data sets with different types of nonlinear relationships to illustrate the methods. All sets consist of $N=10^5$ samples drawn from a certain distribution. The dimension $p$ is allowed to vary from $1$ to $16$. The first set is the independent beta distribution in $p$ dimensions, the second set is a multivariate Gaussian distribution in $p$ dimensions, and the third set is an artificial data set which contains strongly nonlinear relationships between the variables. The datasets are re-scaled to the domain $[0,1]^p$ because the Genz test functions are defined on the unit cube. Their parameters are given as follows.

The beta distribution has parameters $\alpha=2$ and $\beta=5$ and its probability distribution function for one dimension is given by
\begin{equation}
f(x) = \frac{1}{B(\alpha,\beta)}x^{\alpha-1}(1-x)^{\beta-1},
\end{equation}
in which $B(\alpha,\beta)$ is the beta function. In higher dimensions, a tensor product of the $1$-dimensional distribution is used.

The multivariate Gaussian distribution has zero mean, unit variance and correlation coefficients $\sigma_{ij}$ between dimensions $i$ and $j$ given by
\begin{equation}
\sigma_{ij} = \frac{1}{|i-j|+1}.
\end{equation}
This is chosen such that neighboring dimensions have larger correlation coefficients than dimensions far apart.

The third and last distribution is given as
\begin{equation}
\begin{bmatrix}X_1 \\ X_2 \\ \vdots \\ X_p\end{bmatrix} = \begin{bmatrix}U(-2,2) \\ X_1^2 \\ \vdots \\ X_1^p\end{bmatrix} + \sigma N(\mathbf{0},I),
\end{equation}
in which $U(-2,2)$ is the uniform distribution on $[-2,2]$, $\sigma$ is chosen to be $0.5$ and $N(\mathbf{0},I)$ the multivariate standard normal distribution. We refer to this distribution as the ``polynomial distribution''.

\begin{figure}[ht!]
\centering
\subfloat[Beta distributed data]{\includegraphics[width=0.3\textwidth]{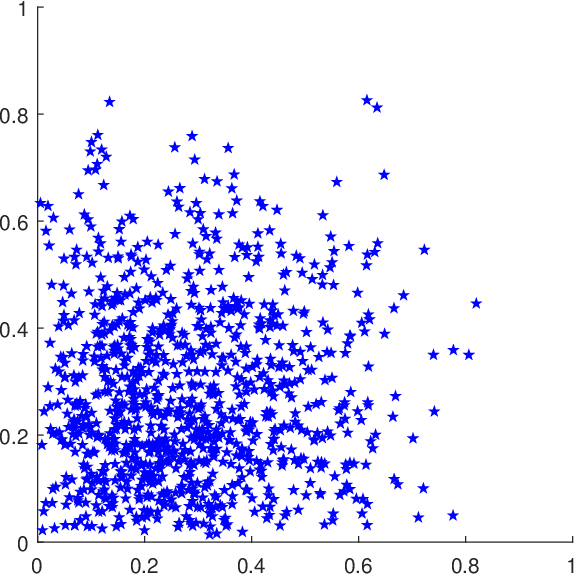}}
\subfloat[Normally distributed data]{\includegraphics[width=0.3\textwidth]{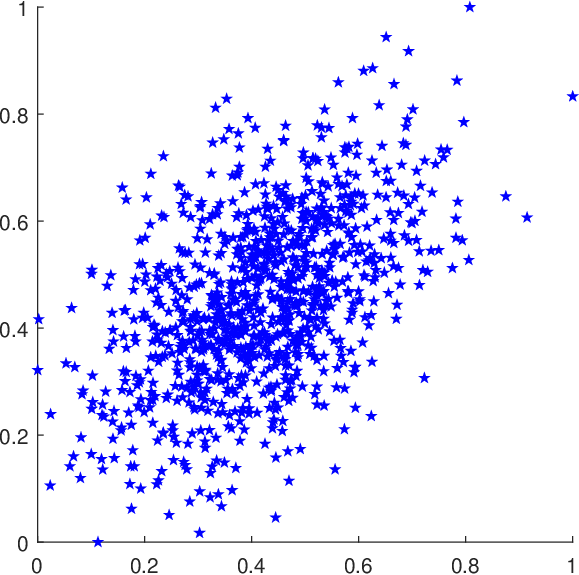}}
\subfloat[Polynomial data]{\includegraphics[width=0.3\textwidth]{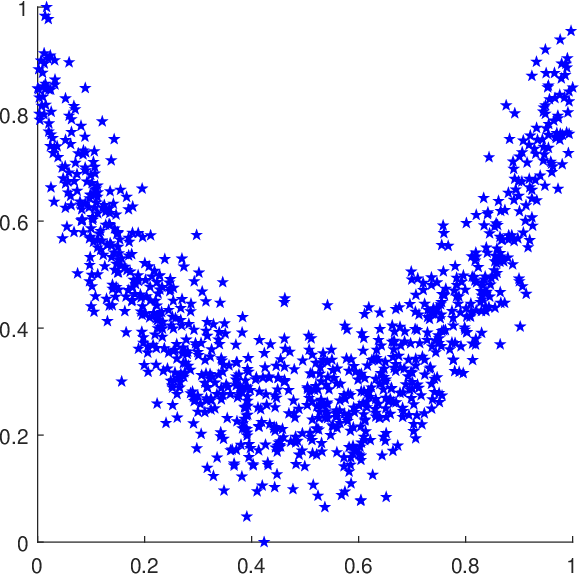}}
\caption{Visualization of the test sets for $p=2$ and $N=10^3$. The beta distributed data is independent, while the normal distributed data is weakly dependent and the polynomial data is strongly, nonlinearly dependent. }
\label{fig:data}
\end{figure}

In Figure \ref{fig:data}, we show $10^3$ data points generated for $p=2$ for the different test sets. From the figure, it is clear that these data sets have different types of nonlinear relationships. The beta distributed data is independent, the normally distributed data is weakly dependent and the polynomial data contains strong nonlinear relationships between the variables and is far from Gaussian.

\begin{figure}[ht!]
\centering
\subfloat[Beta distributed data, $k_{max}=20$]{\includegraphics[width=0.45\textwidth]{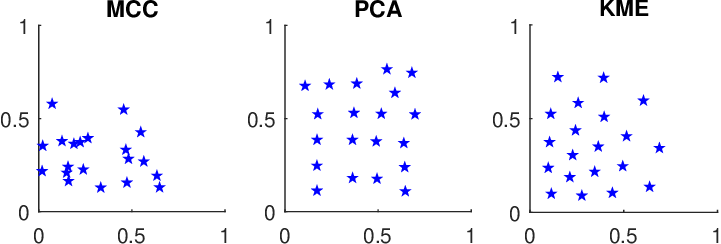}}\hfill
\subfloat[Beta distributed data, $k_{max}=100$]{\includegraphics[width=0.45\textwidth]{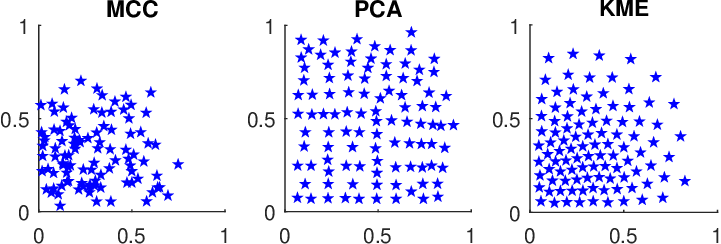}}\\
\subfloat[Normally distributed data, $k_{max}=20$]{\includegraphics[width=0.45\textwidth]{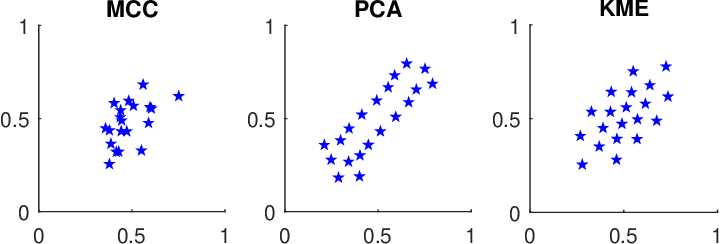}}\hfill
\subfloat[Normally distributed data, $k_{max}=100$]{\includegraphics[width=0.45\textwidth]{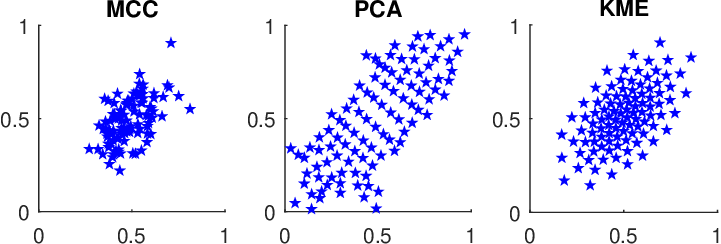}}\\
\subfloat[Polynomial data, $k_{max}=20$]{\includegraphics[width=0.45\textwidth]{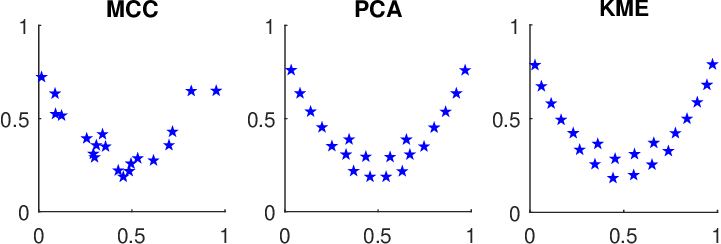}}\hfill
\subfloat[Polynomial data, $k_{max}=100$]{\includegraphics[width=0.45\textwidth]{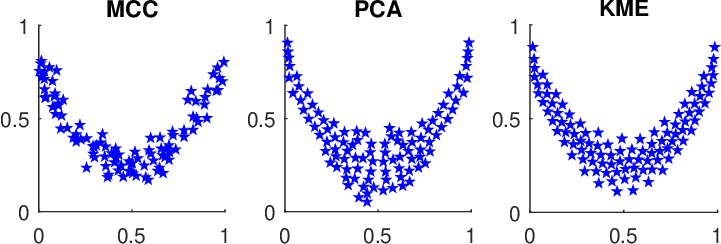}}
\caption{Visualization of the partitionings for $p=2$. The general observation is that MCC and KME have most nodes in dense regions of the data, while PCA is more spread out over the domain of the data. }
\label{fig:cps}
\end{figure}

In Figure \ref{fig:cps}, the partitionings (for $K=20$ and $100$) for the different test sets are shown. One of the observations is that the MCC method yields most clusters in dense regions, just as the KME method. In the latter, the spacing between the nodes is more evenly distributed in space. However, the PCA method also has nodes in less dense regions of the data set and is even more evenly distributed.

\subsection{Tests}
The tests of the methods will consist of integrating the test functions on each of the data sets and comparing the integrals to the Monte Carlo integrals. The data sets will be generated only once and reused. The output of each of the methods is the value of the integral of the test function when performed with the cluster points and weights. Not all results will be shown, but we will show representative examples. The error measure we use is the relative error, defined by
\begin{equation}\label{eq:relerr}
\varepsilon = \frac{|I_{PCA}-I_{MC}|}{|I_{MC}|}.
\end{equation}

\subsubsection{Dimension effects}
\begin{figure}[p]
\centering
\subfloat[Beta distributed data, test function 1]{\includegraphics[width=0.45\textwidth]{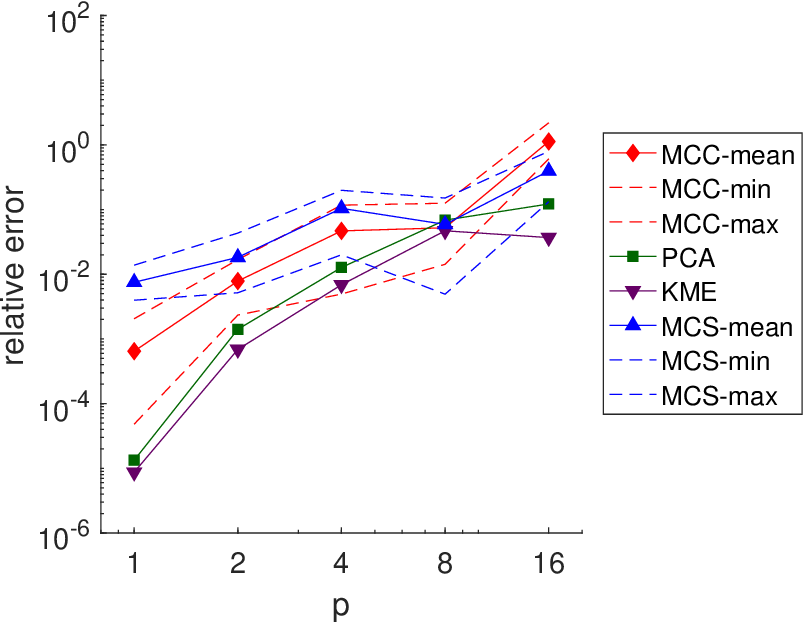}}\hfill
\subfloat[Beta distributed data, test function 2]{\includegraphics[width=0.45\textwidth]{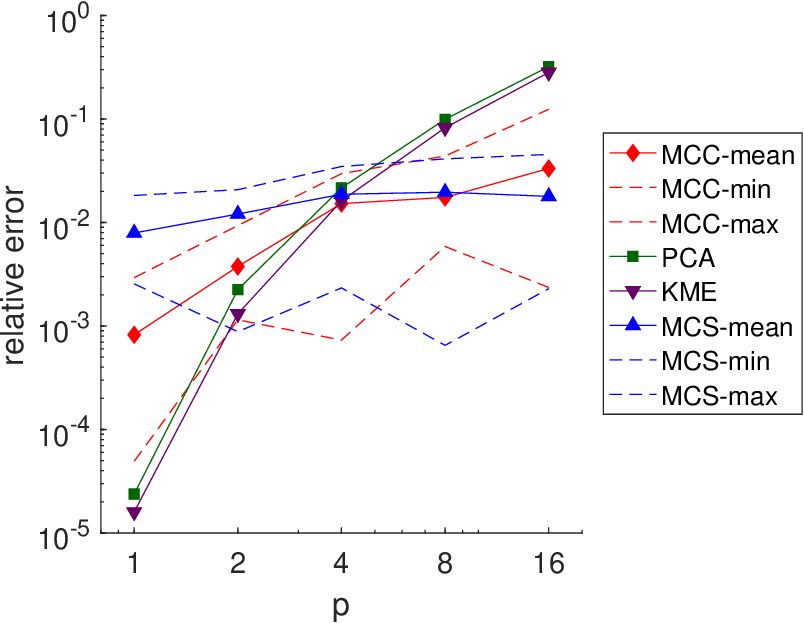}}\\
\subfloat[Normally distributed data, test function 1]{\includegraphics[width=0.45\textwidth]{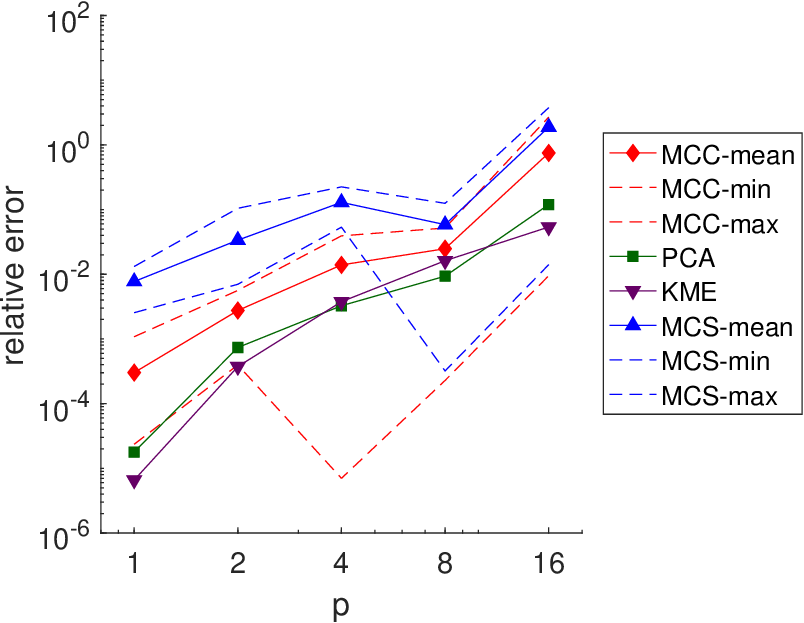}}\hfill
\subfloat[Normally distributed data, test function 2]{\includegraphics[width=0.45\textwidth]{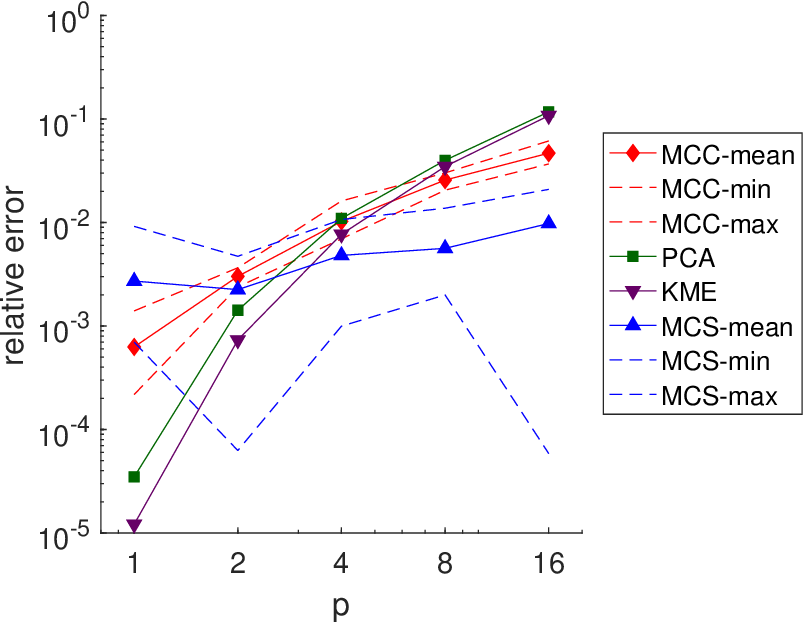}}\\
\subfloat[Polynomial data, test function 1]{\includegraphics[width=0.45\textwidth]{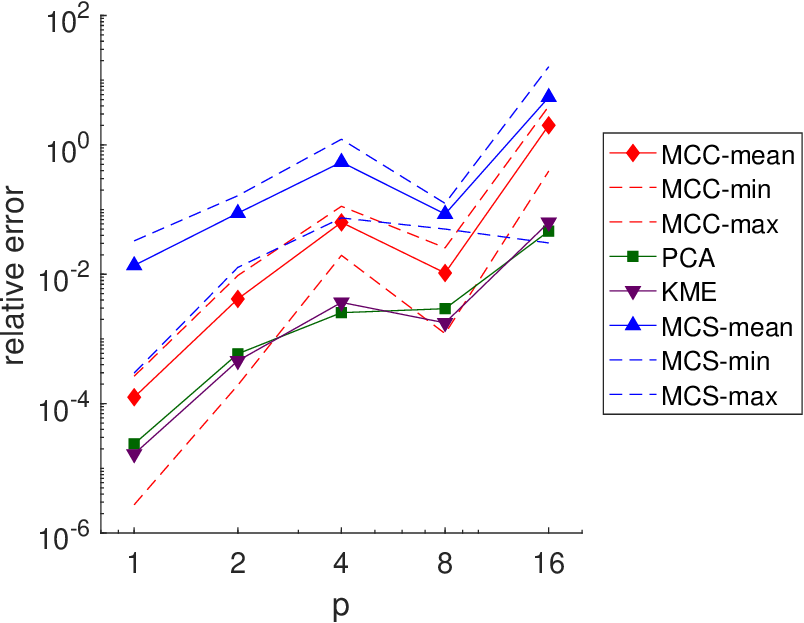}}\hfill
\subfloat[Polynomial data, test function 2]{\includegraphics[width=0.45\textwidth]{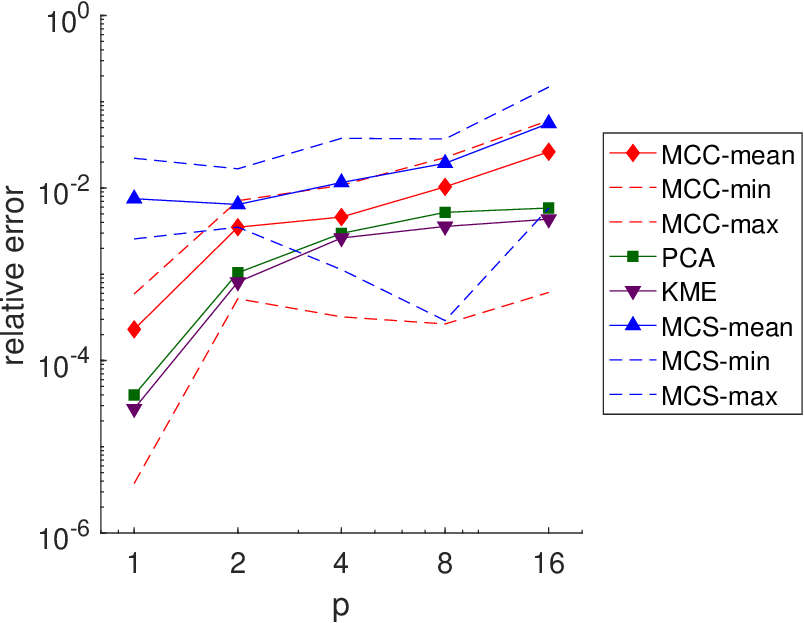}}\\
\caption{Relative error depending on dimension for different methods and data sets with $K=50$.}
\label{fig:relerrfuncp}
\end{figure}

First, we compute the relative error as given by Equation \ref{eq:relerr}. We do this for various values of $p$ and data sets for a fixed maximum number of cluster points $K=50$ to see how the error relates to dimension. The results are in Figure \ref{fig:relerrfuncp} for the first and second test function. In this figure, we observe a trend which holds for all of the proposed methods: namely, that the relative error of PCA and KME are in general lower then for MCC. Furthermore, the error of the MCS does not vary to a large extent with dimension, as expected. Also, it is visible that the PCA-method for the second test function performs better for the dependent data sets and especially for the polynomial data set, which is highly dependent. This indicates that the methods work especially well for dependent inputs, which is caused by the data being concentrated on or near a low-dimensional manifold. It can also be seen that for the deterministic PCA method, the result is more robust with respect to increase of the dimension.
For the test functions $3-5$, the results are similar (not shown). For the discontinuous test function $6$, results are less robust (not shown), due to the discontinuity of the test function. It can also be seen that in some cases, the MCC and MCS results are better than the PCA and KME results, however the variance of the error with MCC and MCS can be large.

\subsubsection{Effect of number of clusters}\label{ssec:k}
\begin{figure}[p]
\centering
\subfloat[Beta distributed data, test function 1]{\includegraphics[width=0.45\textwidth]{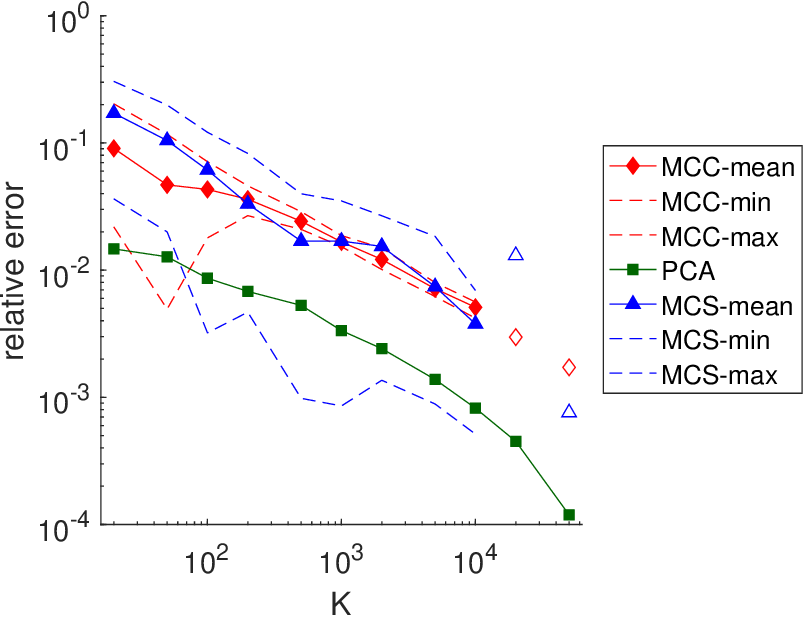}}\hfill
\subfloat[Beta distributed data, test function 2]{\includegraphics[width=0.45\textwidth]{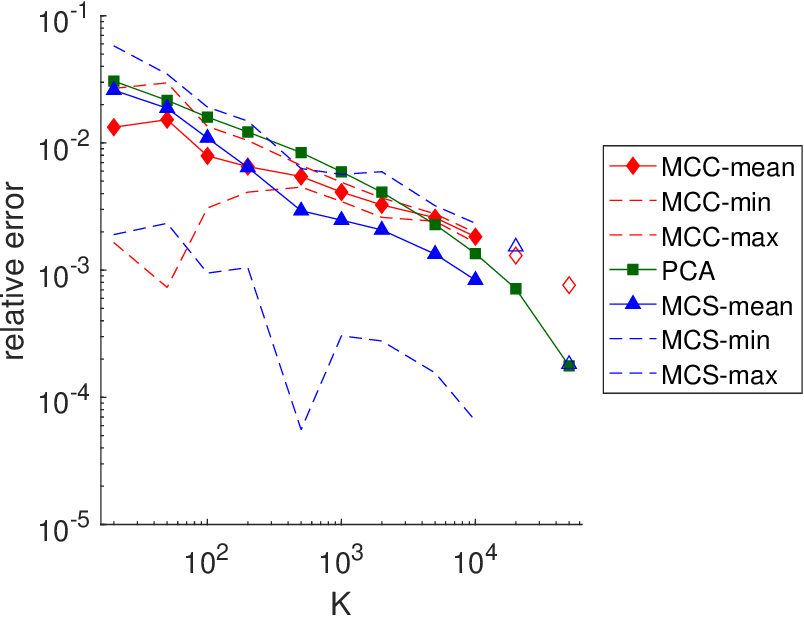}}\\
\subfloat[Normally distributed data, test function 1]{\includegraphics[width=0.45\textwidth]{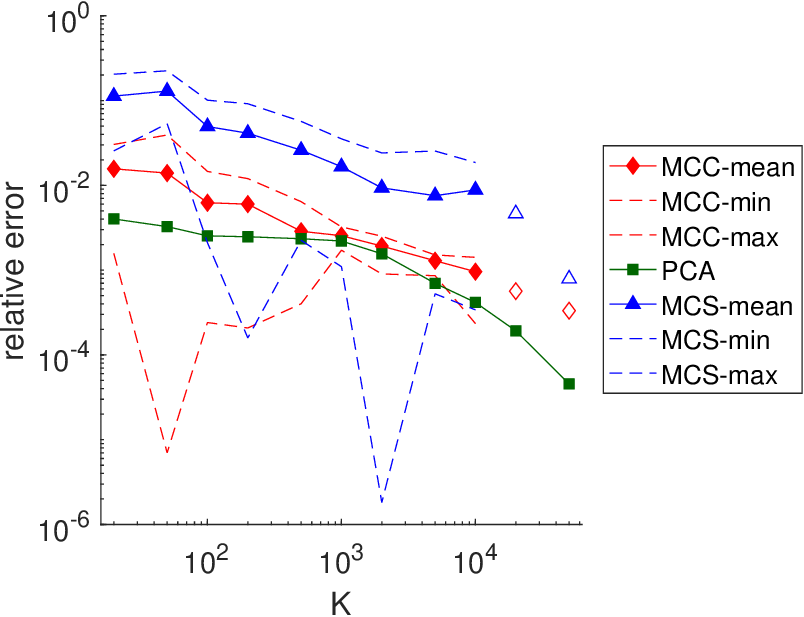}}\hfill
\subfloat[Normally distributed data, test function 2]{\includegraphics[width=0.45\textwidth]{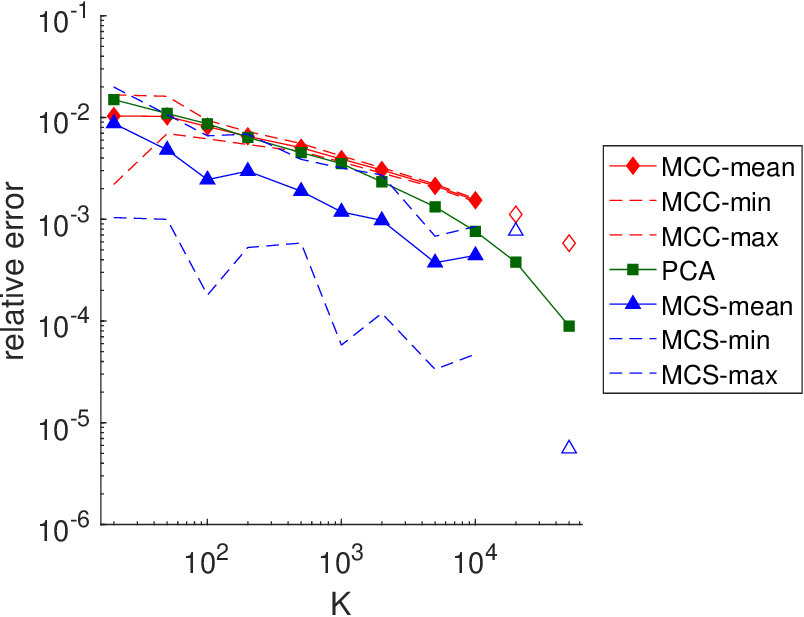}}\\
\subfloat[Polynomial data, test function 1]{\includegraphics[width=0.45\textwidth]{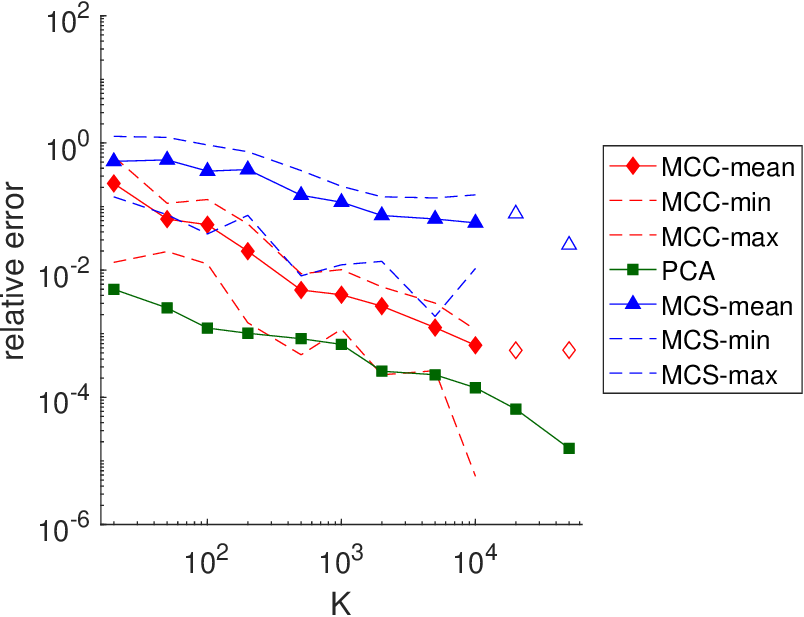}}\hfill
\subfloat[Polynomial data, test function 2]{\includegraphics[width=0.45\textwidth]{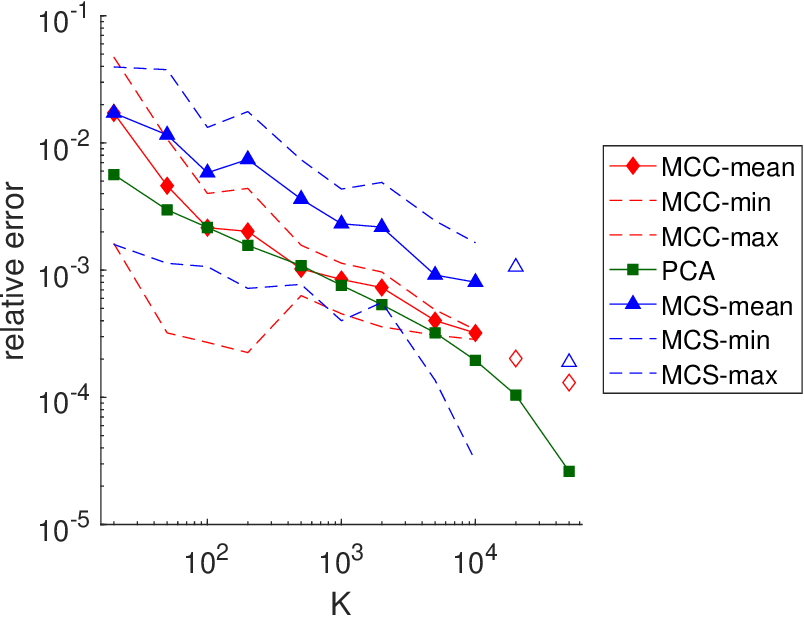}}\\
\caption{Relative error depending on $k_{max}$ for MCC, PCA and MCS and the three data sets with $p=4$. }
\label{fig:relerrfunck}
\end{figure}

In Figure \ref{fig:relerrfunck}, the effect of increasing $K$ is studied for the MCC and PCA-method for $p=4$. These results support the statements from Section \ref{ssec:convergencetheory}, namely that the errors generally decrease with increasing $K$. We show the results for the first and second test function. It can be seen that for the first test function, the PCA-method performs clearly better than the MCC method, which in turn performs better than MCS.
Similar to the previous test, we observe that the methods work better for more strongly dependent datasets. For test function $2$ the PCA, MCC and MCS methods are closer in performance, however the error decrease with increasing $K$ is more robust for PCA.
For the last two settings of $K$, the simulation times restricted the number of simulations to only $1$. Therefore, no minimum and maximum are shown, and the marker is adapted.

\subsubsection{Computational cost}
For MCC and PCA, the time to compute $100$ nodes and weights is in the order of seconds, while it is in the order of minutes for KME (when $25$ repetitions are used to compute one set). For MCS, it is negligible. The time to perform the clustering is about linear in $K$ and $N$ for all methods, although the constants differ. PCA is fastest, followed by KME (which depends on the number of repetitions $r$), while MCC is slowest. This is due to the cost from assigning all data points to clusters and the implementation of the methods.
We emphasize that the clustering needs to be carried out only once, as a pre-processing step to determine the nodes (and associated weights) in parameter space at which the expensive model $u(\mathbf{x})$ must be evaluated. For certain applications, a single evaluation of $u(\mathbf{x})$ can take hours of computation, hence the computational cost of a pre-processing step that takes only seconds to minutes is negligible.

\section{Lid-driven cavity flow}\label{sec:ldcf}
\begin{figure}[ht!]
  \centering
  \includegraphics[width=0.4\textwidth]{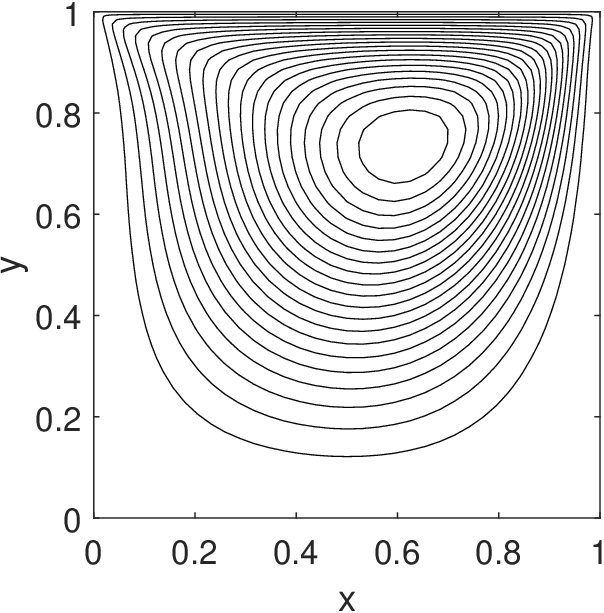}
  \caption{Sketch of the situation (with streamlines). The top wall is moving, while the others contain a no-slip boundary condition. }\label{fig:sketch}
\end{figure}

\begin{figure}[ht!]
  \centering
  \subfloat[Data]{\includegraphics[width=0.32\textwidth]{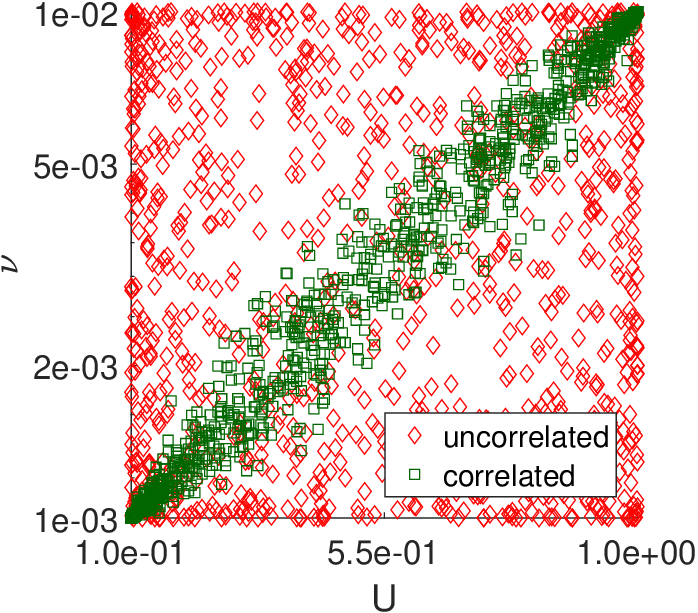}}
  \subfloat[Cluster centers]{\includegraphics[width=0.32\textwidth]{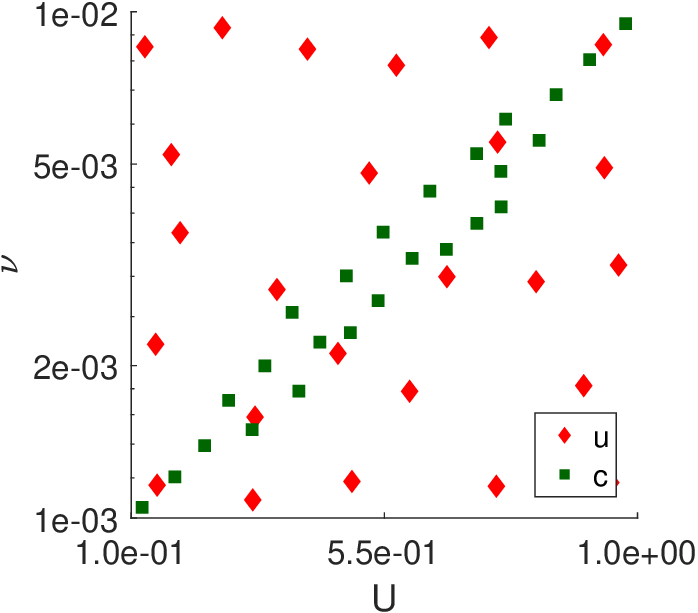}}
  \subfloat[Random samples]{\includegraphics[width=0.32\textwidth]{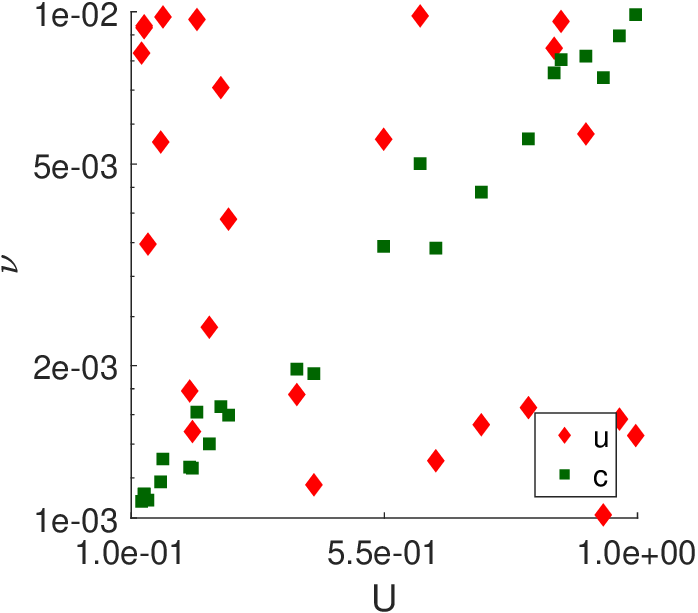}}
  \caption{Visualization of the data,the cluster centers and random samples for the lid-driven cavity flow data.}\label{fig:unu}
\end{figure}

The lid-driven cavity flow is a well known example in computational fluid dynamics for validating new computing methods \cite{2Bur66,Ghi82,2Bot98,2Ert05,2Bru06}. The problem involves fluid flow in a simple, 2D geometry with equally simple boundary conditions. The geometry consists of a (square) box $D=[0,1]^2$ with three fixed walls, and the top wall is moving in one direction with a fixed velocity $U$. The box contains a fluid with viscosity $\nu$ and the incompressible Navier-Stokes equations are solved for the stationary case. The output we consider is the velocity along the centerline at $x=0.5$. The code from \cite{2San11} and \cite{San13} was used for the simulations, with a $50\times 50$ nonuniform grid, which is refined at the boundaries of the domain. A sketch of the situation with Reynolds number $\text{Re}=\frac{UL}{\nu}=100$ ($L=1$) can be found in Figure \ref{fig:sketch}.
We treat $U$ and $\nu$ as uncertain input parameters. Our goal is to demonstrate that we can quantify the uncertainty in the output efficiently by using the PCA method instead of Monte Carlo simulations. To do this, we construct two data sets with $N=10^3$ samples of the velocity $U$ and the viscosity $\nu$. In one data set, they are independent, while in the other set, they are dependent. The samples are generated from a standard Gaussian copula with $\rho=0$ (independent case) and $\rho=-0.99$ (dependent case) and transformed to samples for $U$ and $\nu$ in the following way:
\begin{equation}
  U_i = 0.1+0.9\cdot F_\beta^{-1}(\omega_{i,1}), \nu_i = 10^{-2-F_\beta^{-1}(\omega_{i,2})},
\end{equation}
where $i$ indicates the $i$th sample, $\omega_i=(\omega_{i,1},\omega_{i,2})$ are elements of the Gaussian copula and $F_\beta^{-1}(\cdot)$ is the inverse cumulative distribution function of the beta distribution with parameters $\alpha=\beta=1/2$. This is chosen such that the flow is laminar, different flow profiles occur and the convergence to steady state flow is fast.
On these two data sets, we apply both the PCA based clustering method to get $K=25$ clusters and Monte Carlo sampling to get $K=25$ samples. The sampling is repeated $r=10$ times. The data points, cluster centers and a possible set of samples are in Figure \ref{fig:unu}.

\begin{figure}[ht!]
  \centering
  \includegraphics[width=0.4\textwidth]{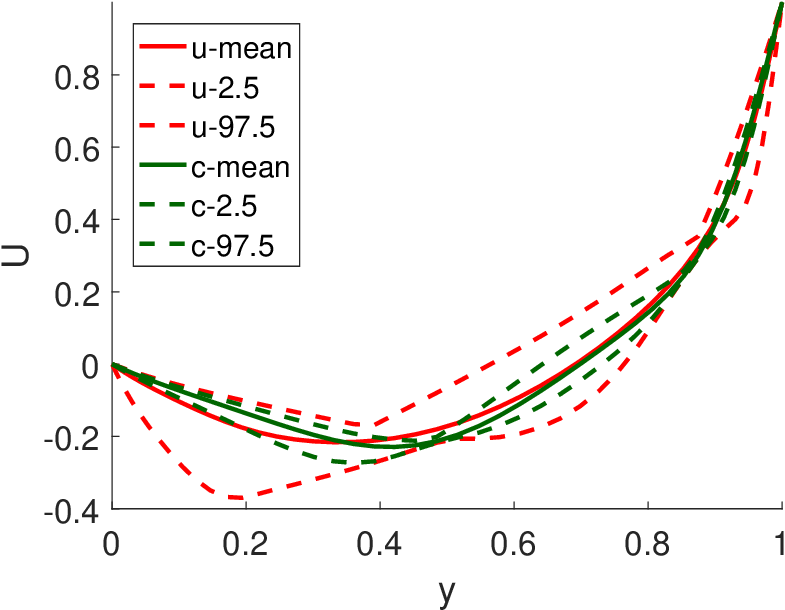}
  \caption{Reference results for the centerline velocity.}\label{fig:diffcu}
\end{figure}

\begin{figure}[ht!]
  \centering
  \subfloat[Independent inputs $(U,\nu)$]{\includegraphics[width=0.45\textwidth]{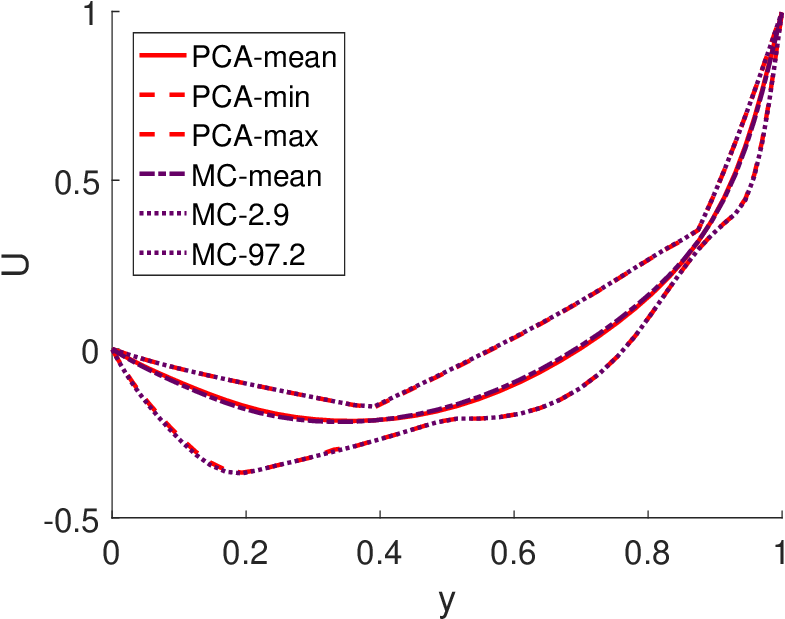}}\hspace{0.05\textwidth}
  \subfloat[Dependent inputs $(U,\nu)$]{\includegraphics[width=0.45\textwidth]{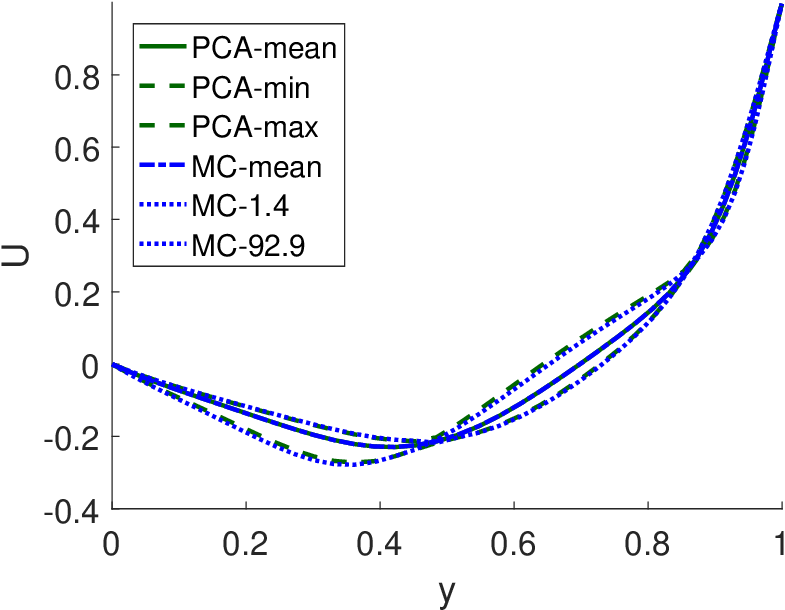}}
  \caption{Centerline velocity computed by the PCA method and compared to the Monte Carlo results.}\label{fig:(un)corPCA}
\end{figure}

\begin{figure}[ht!]
  \centering
  \subfloat[Independent inputs $(U,\nu)$]{\includegraphics[width=0.45\textwidth]{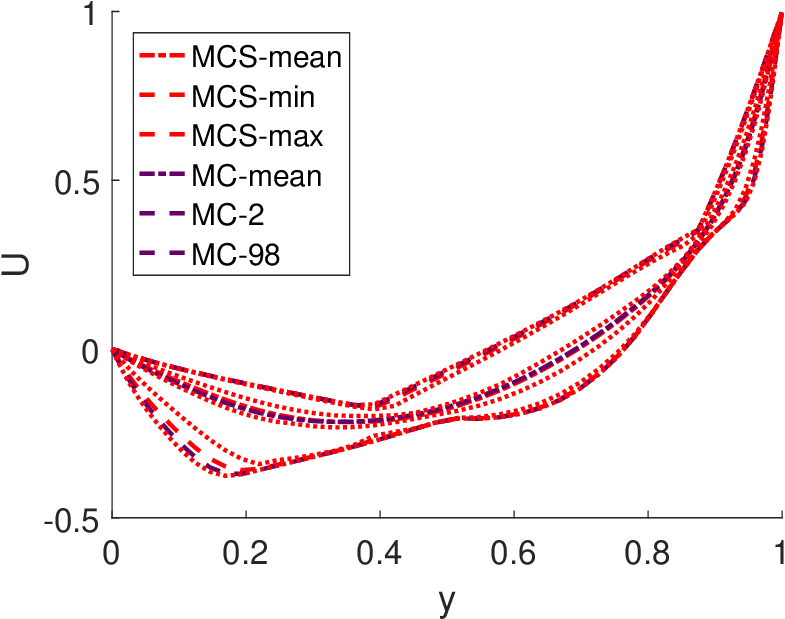}}\hspace{0.05\textwidth}
  \subfloat[Dependent inputs $(U,\nu)$]{\includegraphics[width=0.45\textwidth]{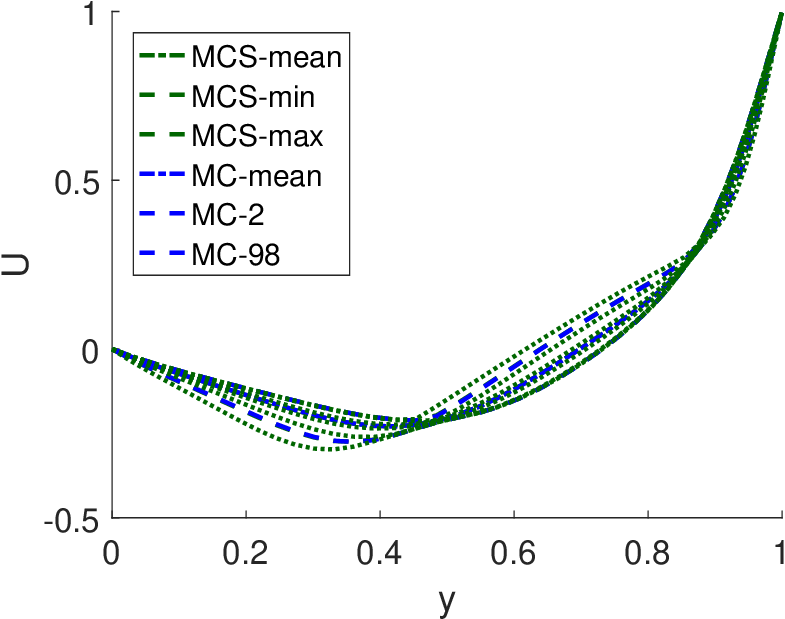}}
  \caption{Centerline velocity computed by Monte Carlo sampling and compared to the Monte Carlo results.}\label{fig:(un)corMC}
\end{figure}

The centerline velocity is computed for the complete data sets and is shown in Figure \ref{fig:diffcu} together with the $2.5$ and $97.5$ percentiles.
We then computed the centerline velocity for the parameter settings given by the $25$ cluster centers and compared the results in Figure \ref{fig:(un)corPCA} for the PCA method.
For the cluster centers, we could not compute the $2.5$ and $97.5$ percentiles exactly due to the cluster weights. Therefore, we show the $2.9$ and $97.2$ percentiles, based on the weight of the cluster centers leading to low or high values of the Reynolds number for the independent data and the $1.4$ and $92.9$ percentiles for the dependent data.
The results in terms of mean values for the statistics mean, minimum and maximum for the Monte Carlo sampling are given in Figure \ref{fig:(un)corMC}.
For these results, also the minimum and maximum computed value of the statistics are given by dotted lines (not in legend), since the experiment was repeated $10$ times. Although the result for the dependent inputs is quite similar, differences exist for the minimum of the independent inputs. Because of the equal weights of all the samples, we used the $2.0$ and $98.0$ percentiles of the data set.

The results match very well, the PCA method required only $25$ evaluations of the fluid flow solver, compared to $10^3$ evaluations used for the full Monte Carlo results.

\section{Conclusion}\label{sec:conclusion}
We have proposed a novel collocation method that employs clustering techniques, thereby successfully dealing with the case of multivariate, dependent inputs. We have assessed the performance of this clustering-based collocation method using the Genz test functions as well as a CFD test case (lid-driven cavity flow) as benchmarks. Three clustering techniques were considered in this context, namely the Monte Carlo (MCC), $k$-means (KME) and principal component analysis based (PCA) clustering techniques. No exact knowledge of the input distribution is needed for the clustering-based method proposed here; a sample of input data is sufficient. Furthermore, for strongly dependent inputs the methods show good performance with input dimension up to $16$. We hypothesize that the more strongly the inputs are dependent, the more the input data are concentrated on a low-dimensional manifold. This makes it possible to obtain a good representation of the input data with a relatively small number of cluster centers.

We observed that the nodes obtained with MCC are mostly concentrated in regions of high density of the input probability distribution, with poor representation of the tails. As a result, this method does not perform well. The PCA method is better at giving a good spread of the collocation nodes, with KME having results in between MCC and PCA. Concerning computational cost, the PCA method is fastest. Overall, the computational cost of the clustering methods is small, and will be negligible compared to the computational cost of expensive model evaluations (involving e.g. CFD solvers).

Altogether, we suggest to use the method based on principal component analysis (PCA) from the ones that we tested. This method is deterministic, it is fast to compute and it yields collocation nodes that are well distributed over the input data set. Also, PCA is better able than KME to include effects from regions of the data with low probability but high impact on the resulting moments. PCA performs well on the tests with Genz functions and has good convergence properties for an increasing number of nodes. Also in the CFD test case of lid-driven cavity flow, PCA performed well.

In this paper we have focused on clustering-based quadrature. However, collocation is also frequently used as an approach for obtaining approximations of output functions through interpolation. We anticipate that the clustering approach we have proposed here will prove useful for interpolation purposes as well. When used for interpolation, the moment estimates might be improved as well, since they are currently based on a function approximation which is piecewise constant on the clusters.

Altogether, the results in this study demonstrate that clustering-based collocation is a feasible and promising approach for UQ with correlated inputs. We intend to develop this approach further in the near future.


\section*{Acknowledgements}

Very sadly, Jeroen Witteveen passed away unexpectedly in the early stages of the research reported here. His presence, inspiration and expertise are greatly missed. Jeroen was one of the initiators of the EUROS project, which includes the current work. This research, as part of the EUROS project, is supported by the Dutch Technology Foundation STW, which is part of the Netherlands Organisation for Scientific Research (NWO), and which is partly funded by the Ministry of Economic Affairs.

 \appendix

 \section{Genz test functions}
 \label{app:genz}
 In Table \ref{tab:Genz}, the definitions of the Genz functions are given. The parameters $\mathbf{a}$ can be used to make the function harder or easier to integrate, while $\mathbf{u}$ contains scale parameters. The functions are defined in $p$ dimensions, in which $p\in\mathbb{N}$, on the domain $[0,1]^p$. In all tests, we will choose $a_i=1$ for $i=1,\ldots,p$. We will choose $u_i=1/2$ for $i=1,\ldots,p$ for all functions except for $f_1$,  where we choose $u_1=0$.
 {\renewcommand{\arraystretch}{2}%
\begin{table}[ht!]
\centering
\caption{Definition of the Genz test functions}
\label{tab:Genz}
\begin{tabular}{c|l|rcl}
\hline
Nr & Characteristic & \multicolumn{3}{l}{Function} \\
\hline
1 & Oscillatory & $f_1(\mathbf{x})$ & = & $\cos\left(2\pi u_1+\sum_{i=1}^p a_ix_i\right)$ \\
2 & Gaussian peak & $f_2(\mathbf{x})$ & = & $\exp\left(-\sum_{i=1}^p a_i^2(x_i-u_i)^2\right)$ \\
3 & $C_0$ & $f_3(\mathbf{x})$ & = & $\exp\left(-\sum_{i=1}^pa_i|x_i-u_i|\right)$ \\
4 & Product peak & $f_4(\mathbf{x})$ & = & $\prod_{i=1}^p \left(a_i^{-2}+(x_i-u_i)^2\right)^{-1}$\\
5 & Corner peak & $f_5(\mathbf{x})$ & = & $\left(1+\sum_{i=1}^p a_ix_i\right)^{-p+1}$\\
6 & Discontinuous & $f_6(\mathbf{x})$ & = & $\left\{\begin{array}{llcc} 0 & x_1>u_1 & \text{ or } & x_2>u_2 \\ \exp\left(\sum_{i=1}^p a_ix_i\right) & \text{else} \end{array}\right.$ \\
\end{tabular}
\end{table}\\
{\renewcommand{\arraystretch}{1}%

In Figure \ref{fig:gtf}, the values for the Genz functions on the domain $[0,1]^2$ are visualized. Test function $2$ and $4$ look the same, but are different.
\begin{figure}[ht!]
\centering
\subfloat[Test function 1]{\includegraphics[width=0.45\textwidth]{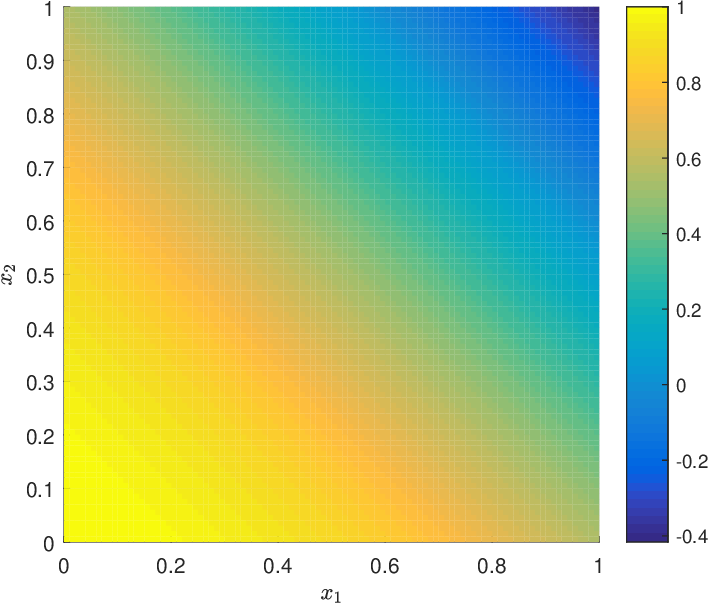}}\hfill
\subfloat[Test function 2]{\includegraphics[width=0.45\textwidth]{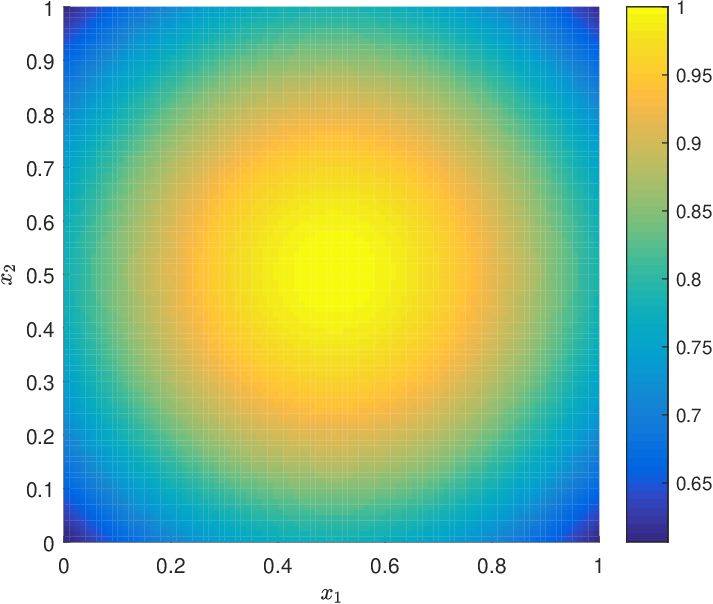}}\\
\subfloat[Test function 3]{\includegraphics[width=0.45\textwidth]{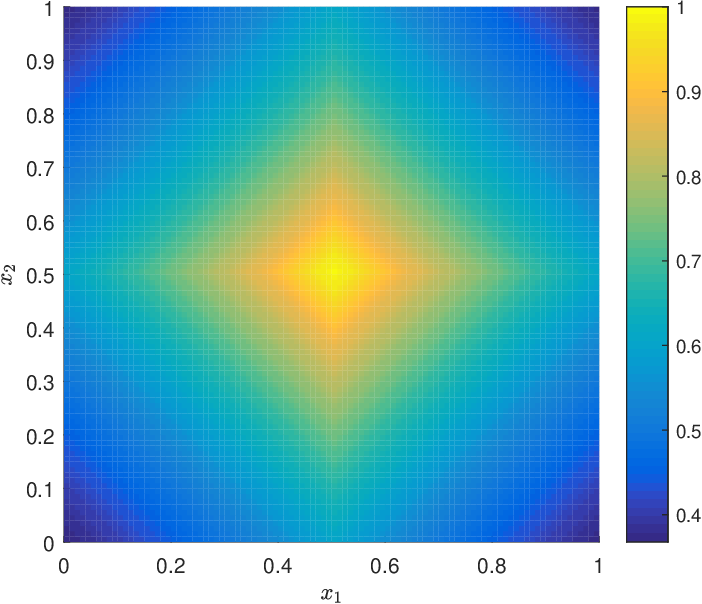}}\hfill
\subfloat[Test function 4]{\includegraphics[width=0.45\textwidth]{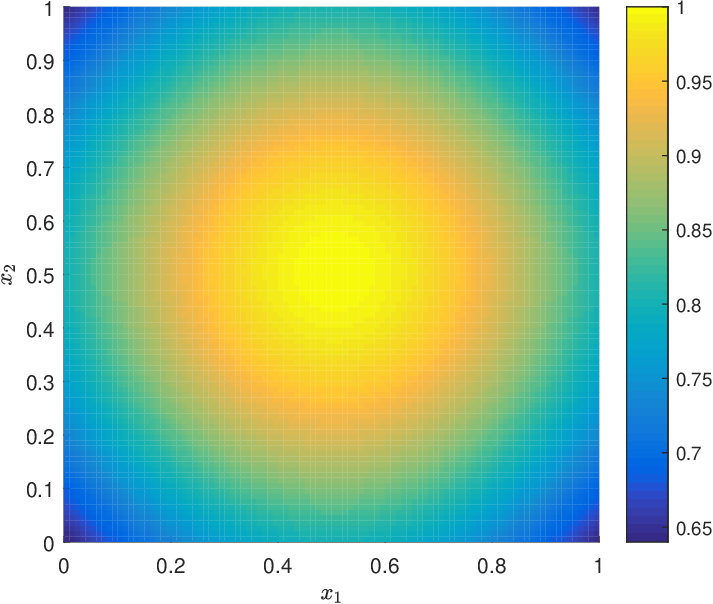}}\\
\subfloat[Test function 5]{\includegraphics[width=0.45\textwidth]{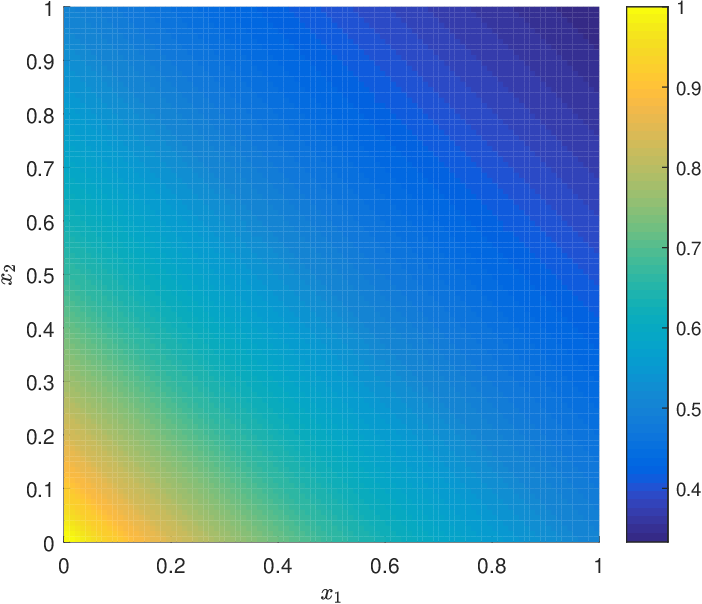}}\hfill
\subfloat[Test function 6]{\includegraphics[width=0.45\textwidth]{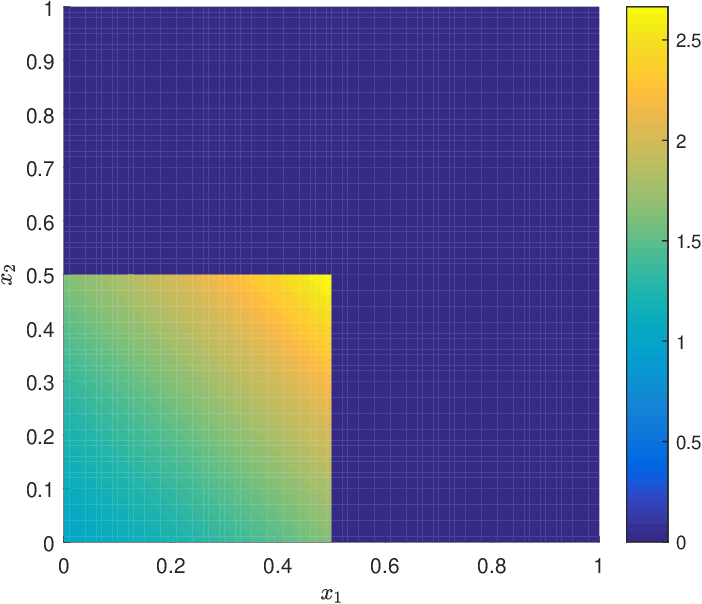}}\\
\caption{Genz test functions.}
\label{fig:gtf}
\end{figure}


\bibliographystyle{IJ4UQ_Bibliography_Style}

\bibliography{references}

\begin{thebibliography}{10}

\bibitem{Bij13}
Bijl, H., Lucor, D., Mishra, S., and Schwab, C., {\em Uncertainty
  Quantification in Computational Fluid Dynamics}, Vol.~92 of Lecture Notes in
  Computational Science and Engineering, Springer, 2013.

\bibitem{Wal02}
Walters, R.W. and Huyse, L., Uncertainty analysis for fluid mechanics with
  applications, Tech.~Rep. 2002-1, ICASE, 2002.

\bibitem{Wit08}
Witteveen, J.A.S. and Bijl, H., Efficient quantification of the effect of
  uncertainties in advection-diffusion problems using polynomial chaos, {\em
  Numerical Heat Transfer, Part B: Fundamentals}, 53(5):437--465, 2008.

\bibitem{Wit07}
Witteveen, J.A.S., Sarkar, S., and Bijl, H., Modeling physical uncertainties in
  dynamic stall induced fluid--structure interaction of turbine blades using
  arbitrary polynomial chaos, {\em Computers \& Structures}, 85(11):866--878,
  2007.

\bibitem{Yil15}
Yildirim, B. and Karniadakis, G.E., Stochastic simulations of ocean waves: An
  uncertainty quantification study, {\em Ocean Modelling}, 86:15--35, 2015.

\bibitem{Xiu02}
Xiu, D. and Karniadakis, G.E., The {W}iener--{A}skey polynomial chaos for
  stochastic differential equations, {\em SIAM {J}ournal on {S}cientific
  {C}omputing}, 24(2):619--644, 2002.

\bibitem{Xiu05}
Xiu, D. and Hesthaven, J.S., High-order collocation methods for differential
  equations with random inputs, {\em SIAM Journal on Scientific Computing},
  27(3):1118--1139, 2005.

\bibitem{Gha03}
Ghanem, R.G. and Spanos, P.D., {\em Stochastic {F}inite {E}lements: {A}
  {S}pectral {A}pproach}, Dover, 2003.

\bibitem{Mai10}
Le~Ma\^{i}tre, O.P. and Knio, O.M., {\em Spectral {M}ethods for {U}ncertainty
  {Q}uantification: {W}ith {A}pplications to {C}omputational {F}luid
  {D}ynamics}, Scientific {C}omputation, Springer, 2010.

\bibitem{Eld09}
Eldred, M.S. and Burkardt, J., Comparison of non-intrusive polynomial chaos and
  stochastic collocation methods for uncertainty quantification, {\em AIAA
  paper 2009-976}, 2009.

\bibitem{2Ros52}
Rosenblatt, M., Remarks on a multivariate transformation, {\em The annals of
  mathematical statistics}, 23(3):470--472, 1952.

\bibitem{Nav15}
Navarro, M., Witteveen, J.A.S., and Blom, J., Stochastic collocation for
  correlated inputs, In {\em UNCECOMP 2015}, 2015.

\bibitem{Smo63}
Smolyak, S.A., Quadrature and interpolation formulas for tensor products of
  certain classes of functions, {\em Sov. Math. Dokl.}, 4:240--243, 1963.

\bibitem{2Ger98}
Gerstner, T. and Griebel, M., Numerical integration using sparse grids, {\em
  Numerical algorithms}, 18(3-4):209--232, 1998.

\bibitem{Ste56}
Steinhaus, H., Sur la division des corps mat\'{e}riels en parties, {\em
  Bulletin de l'Acad\'{e}mie {P}olonaise des {S}ciences}, IV(12):801--804,
  1956.

\bibitem{Jai99}
Jain, A.K., Murty, M.N., and Flynn, P.J., Data clustering: a review, {\em ACM
  computing surveys (CSUR)}, 31(3):264--323, 1999.

\bibitem{Mac67}
MacQueen, J. , Some methods for classification and analysis of multivariate
  observations, In {\em Proceedings of the {F}ifth Berkeley {S}ymposium on
  {M}athematical {S}tatistics and {P}robability}, Vol.~1, pp. 281--297, 1967.

\bibitem{Elk03}
Elkan, C., Using the triangle inequality to accelerate $k$-means, In {\em
  Proceedings ICML-2003}, Vol.~3, pp. 147--153, 2003.

\bibitem{Lik03}
Likas, A., Vlassis, N., and Verbeek, J.J., The global $k$-means clustering
  algorithm, {\em Pattern {R}ecognition}, 36(2):451--461, 2003.

\bibitem{Bag08}
Bagirov, A.M., Modified global $k$-means algorithm for minimum sum-of-squares
  clustering problems, {\em Pattern Recognition}, 41(10):3192--3199, 2008.

\bibitem{Han05}
Hansen, P., Ngai, E., Cheung, B.K., and Mladenovic, N., Analysis of global
  $k$-means, an incremental heuristic for minimum sum-of-squares clustering,
  {\em Journal of {C}lassification}, 22(2):287--310, 2005.

\bibitem{Coh14}
Cohen, M.B., Elder, S., Musco, C., Musco, C., and Persu, M., Dimensionality
  reduction for $k$-means clustering and low rank approximation, {\em arXiv
  preprint arXiv:1410.6801}, 2015.

\bibitem{Art07}
Arthur, D. and Vassilvitskii, S., ${K}$-means++: The advantages of careful
  seeding, In {\em Proceedings of the Eighteenth Annual ACM-SIAM Symposium on
  Discrete Algorithms}, pp. 1027--1035. SIAM, 2007.

\bibitem{Din04}
Ding, C. and He, X., {$K$}-means clustering via principal component analysis,
  In {\em Proceedings of the twenty-first international conference on Machine
  learning}, p.~29. ACM, 2004.

\bibitem{Su04}
Su, T. and Dy, J., A deterministic method for initializing k-means clustering,
  In {\em Proceedings ICTAI 2004}, pp. 784--786. IEEE, 2004.

\bibitem{Gue91}
Gu{\'e}noche, A., Hansen, P., and Jaumard, B., Efficient algorithms for
  divisive hierarchical clustering with the diameter criterion, {\em Journal of
  {C}lassification}, 8(1):5--30, 1991.

\bibitem{Fan08}
Fang, H. and Saad, Y., Farthest centroids divisive clustering, In {\em
  Proceedings ICMLA'08}, pp. 232--238. IEEE, 2008.

\bibitem{Gen84}
Genz, A., Testing multidimensional integration routines, In {\em Proceedings of
  International Conference on Tools, Methods and Languages for Scientific and
  Engineering Computation}, pp. 81--94. Elsevier North-Holland, 1984.

\bibitem{2Bur66}
Burggraf, O.R., Analytical and numerical studies of the structure of steady
  separated flows, {\em Journal of {F}luid {M}echanics}, 24(1):113--151, 1966.

\bibitem{Ghi82}
Ghia, U.K.N.G., Ghia, K.N., and Shin, C.T., High-{R}e solutions for
  {I}ncompressible {F}low {U}sing the {N}avier-{S}tokes {E}quations and a
  {M}ultigrid {M}ethod, {\em Journal of {C}omputational {P}hysics},
  48(3):387--411, 1982.

\bibitem{2Bot98}
Botella, O. and Peyret, R., Benchmark spectral results on the lid-driven cavity
  flow, {\em Computers \& {F}luids}, 27(4):421--433, 1998.

\bibitem{2Ert05}
Erturk, E., Corke, T.C., and G{\"o}k{\c{c}}{\"o}l, C., Numerical solutions of
  2-{D} steady incompressible driven cavity flow at high {R}eynolds numbers,
  {\em International {J}ournal for {N}umerical {M}ethods in {F}luids},
  48(7):747--774, 2005.

\bibitem{2Bru06}
Bruneau, C.H. and Saad, M., The 2{D} lid-driven cavity problem revisited, {\em
  Computers \& F{}luids}, 35(3):326--348, 2006.

\bibitem{2San11}
Sanderse, B., {E}nergy-{C}onserving {N}avier-{S}tokes {S}olver. {V}erification
  of steady laminar flows, Tech.~Rep. E-11-042, ECN, 2011.

\bibitem{San13}
Sanderse, B., Energy-conserving discretization methods for the incompressible
  {N}avier-{S}tokes equations, PhD thesis, Eindhoven University of Technology,
  2013.

\end{thebibliography}
\end{document}